\documentclass[11pt,a4paper]{article}
\usepackage{amsmath,amssymb,latexsym,amsthm}
\usepackage[english]{babel}
\usepackage{graphicx}
\usepackage{hhline}
\newtheorem{Theorem}{Theorem}
\theoremstyle{remark}
\newtheorem{Remark}[Theorem]{Remark}
\newtheorem{Example}[Theorem]{Example}
\newcommand{\Pro}{\noindent{\em Proof. }}

\begin{document}
\title{\vspace{1cm} The action of a compact Lie group on nilpotent Lie algebras of type $\{n,2\}$
\thanks{2010 Mathematics Subject Classification: 17B05, 17B30, 15A63}
\thanks{Keywords: Nilpotent Lie algebras, Pair of alternating forms.}
\thanks{This research was supported by Universit\`a di Palermo (2012-ATE-0446) and by the European Union's Seventh Framework Programme (FP7/2007-2013) under grant agreements no. 317721, no. 318202. }
}
\author{Giovanni Falcone and \'Agota Figula}
\date{}
\maketitle 

\begin{abstract}
\noindent
We classify  finite-dimensional nilpotent Lie algebras with  $2$-dimensional central commutator ideals admitting a Lie group of automorphisms isomorphic to ${\mathrm{SO}}_2(\mathbb R)$. This enables one to enlarge the class of
nilpotent Lie algebras of type $\{n,2\}$.
\end{abstract}

\centerline{\bf 1. Introduction}

\bigskip
The most simple (non-abelian) Lie algebras are the generalized Heisenberg Lie algebras, defined on a $(n+1)$-dimensional vector space
${\mathfrak h}=V\oplus\langle x\rangle$ by a non-degenerate alternating form $F$ on the $n$-dimensional subspace $V$ ($n$ even),
putting $[u,v]=F(u,v)x$, for any $u,v\in V$.

According to the literature beginning with Vergne \cite{Vergne}, metabelian  Lie algebras ${\mathfrak h}=V\oplus\langle x,y\rangle$ of dimension $(n+2)$ defined by a pair of
alternating forms $F_1, F_2$ on the $n$-dimensional vector space $V$, putting, for any $u,v\in V$, 
$[u,v]=F_1(u,v)x+F_2(u,v)y$, are called  nilpotent Lie algebras of type $\{n,2\}$, where the \emph{type} 
$\{p_1,\dots, p_c\}$ of a nilpotent Lie algebra ${\mathfrak g}$ with descending central series 
${\mathfrak g}^{(i)}=[{\mathfrak g},{\mathfrak g}^{(i-1)}]$ is defined by the integers
$p_i={\mathrm{dim}}{\frac{{\mathfrak g}^{(i-1)}}{{\mathfrak g}^{(i)}}}$.

Nilpotent (real or complex) Lie algebras of type $\{n,2\}$ have been classified firstly by Gauger \cite{Gauger}, applying the canonical reduction of the pair $F_1,F_2$. We mention that also nilpotent Lie algebras of type $\{n,2, 1\}$ can be explicitly described (cf. 
\cite{Bartolone}).
According to results of Belitskii, Lipyanski, and Sergeichuk \cite{Belitskii}, this line of investigation cannot be carried further. A possible way of broadening these families of Lie algebras appears therefore that of considering their derivations.

In this paper we want to study derivations of a nilpotent Lie algebra ${\mathfrak h}$ of type $\{n,2\}$, whereas derivations of a nilpotent Lie algebra of type $\{n,2, 1\}$ are being considered in another paper \cite{Bartolone2}.

As a first question, we ask whether ${\mathfrak h}$ admits a compact Lie algebra of derivations. This question is interesting for the study of the isometry groups of homogeneous nilmanifolds.    
Considering an invariant inner product on ${\mathfrak h}$ the compact Lie algebra of derivations of ${\mathfrak h}$  belongs to the Lie algebra of the isometry group of a simply connected nilmanifold associated to ${\mathfrak h}$ 
(cf. \cite{wilson}, p. 337).  
Since a non-commutative simple compact Lie algebra cannot have a two-dimensional representation, a non-commutative simple compact Lie algebra of derivations of a nilpotent Lie algebra ${\mathfrak h}$ of type $\{n,2\}$ must induce the null map on the two-dimensional commutator ideal ${\mathfrak h}'$. 

The smallest example, that is, the $5$-dimensional Lie algebra of type $\{3,2\}$ defined by
$$[u_1,u_2]=x, \quad [u_1,u_3]=y,$$
does not have compact non-commutative Lie algebras of derivations, since its derivations inducing the null map on  
${\mathfrak h}'$ are defined with respect to the basis $\{u_1,u_2,u_3,x,y\}$  by the matrices
$$\left(\begin{array}{c|cc|cc}a&0&0&0&0\\
\hline b&-a&0&0&0\\
c&0&-a&0&0\\
\hline d_1&d_2&d_3&0&0\\
d_4&d_5&d_6&0&0\end{array}\right).$$
But its group of automorphisms contains the group $SO_2(\mathbb R)$.

In general the structure of a Lie algebra ${\mathfrak h}$ of type $\{n,2\}$ is not particularly rigid, and
the following example shows that, as soon as $n=4$, the algebra of derivations contains compact simple subalgebras.

\begin{Example}\label{example}
Let  ${\mathcal {B}}=\{u_1,u_2,u_3,u_4,x,y\}$ be a basis of the $6$-dimensional Lie algebra $\mathfrak h$ of type 
$\{4,2\}$ defined by
$$[u_1,u_3]=x, \quad [u_1,u_4]=-y,\quad [u_2,u_3]=y,\quad [u_2,u_4]=x.$$
 A direct computation shows that, with respect to the basis ${\mathcal {B}}$, the derivations of $\mathfrak h$ inducing the null map on
 ${\mathfrak h}'$ are represented by matrices of the form
$$\left(\begin{array}{cccc|cc}{ a_1}&{ a_2}&{ a_3}&{ a_4}&0&0\\
-{ a_2}& { a_1}&{ a_4}&-{ a_3}&0&0\\
 -{ b_2}&{ c_2}&-{ a_1}&{ a_2}&0&0\\
 { c_2}&{ b_2}&-{ a_2}&-{ a_1} &0&0\\
\hline d_1&d_2&d_3&d_4&0&0\\
d_5&d_6&d_7&d_8&0&0\end{array}\right).$$
With $a_1=0$, $c_2=-a_4$, $b_2=a_3$, and all the entries $d_i$ equal to zero, we get an algebra 
 isomorphic to a real form of ${\mathfrak{su}}_2({\mathbb{C}})$. 
This Lie algebra $\mathfrak h$ is isomorphic to the Lie algebra of the complex Heisenberg group 
$$N=\left\{\left(\begin{array}{ccc} 1 & \alpha & \gamma \\
0 & 1 & \beta \\
0 & 0 & 1 \end{array}\right), \alpha , \beta , \gamma \in \mathbb C \right\}, $$
which is an interesting example of nilmanifolds. Namely, the group $N$ is the smallest $H$-type group. The $H$-type groups with two-dimensional centre are the complex Heisenberg groups (cf. \cite{tamaru}, p. 3252).  
\end{Example}

\smallskip
This, and the fact that any maximal compact subgroup of a connected real solvable Lie group is a torus, legitimates one, in our opinion, to study the action of a torus on nilpotent Lie algebras of type $\{n,2\}$.

The simple structure of a nilpotent Lie algebra of type $\{n,2\}$ admitting a one-dimensional compact group $T$ of automorphisms is hidden by three obstacles that can be removed by a clear notation. The first is the representation of 
$2h\times 2k$ real matrices as $h\times k$ matrices with coefficients in the algebra of split-quaternions. The second is the reduction to canonical form of a pair of alternating forms. The third is the reduction to the $T$-\emph{undecomposable} case. In Section 2 we summarize some known facts, fix the notation and find the relations (\ref{T1T2blocksbeta}), that are basic for the classification. Thereafter, we consider, in Section 3, the case where $T$ induces the identity on the commutator ideal. It turns out that the classification is parametrized only by the dimension of $\mathfrak h$ and by a complex eigenvalue $q$ (cf. Theorem \ref{Thmbeta=0}). In Section 4, we consider the case where $T$ operates effectively on the commutator ideal. In this case, the classification is more rich and we consider four cases. The classification is always reached  by using parameters that are  linked to the dimension of the eigenspaces of $T$ and by a certain arbitrarity in the reduction to the echelon form of the blocks of the matrices describing the Lie algebra $\mathfrak h$. But in one case (cf. Theorem \ref{skewHermitian}), a class is possibly given by the real form of an arbitrary skew-symmetric complex matrix. Our classes give remarkable examples for nilmanifolds $M$ such that the group of isometries of 
$M$ contains the compact group $SO_2(\mathbb R)$.

\bigskip
\noindent
\centerline{\bf 2. Notation} \label{Notation}

\bigskip
\noindent
\textbf{{Split-quaternions and matrix notation.}}
We denote: \begin{itemize}
\item[i)] by $0$ any $n\times m$ (real or complex) zero matrix, by $I_m$ the (real or complex) $m$-dimensional identity matrix, and by $\widetilde{I}_{n\times m}$ a $n\times m$ matrix having rank $m$, obtained by the identity matrix $I_m$ by inserting $n-m$ zero rows (without specifying, however, which ones);
\item[ii)] by  $J=\left(\begin{array}{cc}
0&1\\
-1&0\end{array}\right)$ the real matrix corresponding to the imaginary unit;
\item[iii)] by $A'$ the transpose of $A$, and by $A^\dag$ the conjugate transpose of $A$;
\item[iv)] by $A\oplus B$ the diagonal block matrix $\left(\begin{array}{cc}
A&0\\
0&B\end{array}\right)$, and by $(\oplus A)$ the diagonal block matrix $A\oplus\dots\oplus A$.
\end{itemize}
Throughout the paper, we represent the Clifford algebra of split-quaternions as the set
$${\mathbb{H}}_-=\{z_1+z_2\omega:z_i\in {\mathbb{C}}, \omega z=\bar{z}\omega, \omega^2=1\}.$$
We recall that, through the usual identification of the complex number $z=a+ib$ with the real matrix $\left(\begin{array}{cc}
a&b\\
-b&a\end{array}\right)$ and of the reflection $\Omega=\left(\begin{array}{cc}
0&1\\
1&0\end{array}\right)$ with the split-quaternion $\omega$, one obtains an isomorphism of the algebra of real $2\times 2$ matrices with the algebra ${\mathbb{H}}_-$ and, more generally, of the space of ${\mathbb R}^{2n\times 2m}$ matrices with the space of
${\mathbb H}_-^{n\times m}$ matrices.

In more details, with the identification of the split-quaternion matrix $\omega I_n$ with the $2n\times 2n$ reflection $\Omega_{2n}= \Omega\oplus\dots\oplus \Omega$,
any matrix $A \in{\mathbb R}^{2n\times 2m}$ can be written in a unique way as $A=A_1+A_2\Omega_{2m}$, where $A_1$ and $A_2$ are \emph{real forms}
of complex matrices $\widehat{A}_1=(z_{ij})$, $\widehat{A}_2=(u_{ij}) \in \mathbb C^{n \times m}$ such that, for ${\overline{\widehat{A}}}_1=(\bar z_{ij})$ and
${\overline{\widehat{A}}}_2=(\bar u_{ij})$, one has $\omega I_n \widehat{A}_i={\overline{\widehat{A}}}_i\, \omega I_{m}$\, ($i=1,2$).

\bigskip
\noindent
\textbf{{Canonical form of a pair of alternating forms.}}
The problem of the simultaneous reduction to canonical form of a pair of symmetric or alternating forms is classic and, even if the alternating case has been settled only in 1976 by R. Scharlau \cite{Scharlau}, it goes back to a paper of Kronecker \cite{Kronecker}. A summary is given in \cite{Olomuc}, here we point out that an unordered pair of alternating complex bilinear forms can be simultaneously reduced to the direct sum of the canonical forms $(L_t,R_t)^\nabla$ and $(I_t, J_t(q))^\nabla$,
where
$$(A,B)^\nabla=\left(\left(\begin{array}{c|c}0&-A\\ \hline A'&0\end{array}\right),\left(\begin{array}{c|c}0&-B\\ \hline B'&0\end{array}\right)\right);$$
$L_t$ and $R_t$ are $t\times (t+1)$ matrices mapping the $t$-tuple $(x_1,\dots, x_t)$ to $(x_1,\dots, x_t,0)$ and $(0,x_1,\dots, x_t)$, respectively; $I_t$ is the identity $t\times t$ matrix, and $J_t(q)$ is the $t\times t$ Jordan block of eigenvalue $q$.
Over the field of real numbers, together with $(L_t,R_t)^\nabla$ and $(I_t, J_t(q))^\nabla$ ($q\in{\mathbb R}$), one has to consider also the real form of $(\widehat{I}_{t}, \widehat{J}_t(q))^\nabla$ ($q\in{\mathbb C}$).

\medskip
Let $\{e_1,\dots, e_n,x,y\}$ be a basis of a nilpotent Lie algebra with a central commutator ideal ${\mathfrak{h}}'$, such that $\{x,y\}$ is a basis of ${\mathfrak{h}}'$. We define the pair of alternating matrices  $\Big(T_1=(a_{ij}),T_2=(b_{ij})\Big)$ by
$$[e_i,e_j]=a_{ij}x+b_{ij}y.$$
Manifestly, such nilpotent Lie algebras can be thoroughly described by the canonical forms
$(L_t,R_t)^\nabla$, $(I_t, J_t(q))^\nabla$ ($q\in{\mathbb R}$) and the real form of $(\widehat{I}_{t}, \widehat{J}_t(q))^\nabla$ ($q\in{\mathbb C}$).
Notice that the module $(L_0,R_0)^\nabla=\big((0),(0)\big)$ is a component of the pair $(T_1,T_2)$ if and only if the centre of $\mathfrak h$ contains the commutator ideal properly.

\bigskip
\noindent
\textbf{{$T$-undecomposable nilpotent Lie algebra of type $\{n,2\}$.}}
Let $\mathfrak h$ be a nilpotent Lie algebra of type $\{n,2\}$, that is, with a $2$-dimensional commutator ideal ${\mathfrak{h}}'$ coinciding with the centre $\mathfrak z$.
Let $T$ be a group of automorphisms of $\mathfrak h$ isomorphic to ${\mathrm{SO}}_2(\mathbb R)$ and ${\mathfrak t}$ the corresponding compact algebra of derivations of ${\mathfrak h}$, let $n=2m$ if $n$ is even, and $n=2m+1$ if $n$ is odd, with
$n\ge 3$.
By the complete reducibility of $T$, we find a basis $\{e_1,\dots, e_n,x,y\}$ of $\mathfrak h$, such that $\{x,y\}$ is a basis of ${\mathfrak{h}}'=\mathfrak z$ and such that ${\mathfrak t}$ operates on ${\mathfrak h}$ as the algebra of matrices
\begin{equation}\label{Teven}
\begin{array}{c}
\left\{\partial(t):=\Big(\alpha_1 t\cdot J \oplus\cdots\oplus\alpha_m t\cdot J\Big) \oplus\beta t\cdot J |\; t \in{\mathbb R}\right\}  \hbox{ for } \ n=2m, \\
\left\{\partial(t):= \Big(0 \oplus \alpha_1 t\cdot J \oplus\cdots\oplus\alpha_m t\cdot J\Big) \oplus\beta t\cdot J |\; t \in {\mathbb R}\right\} \  \hbox{for} \ n=2m+1,   \end{array}\end{equation}
where $\beta t\cdot J$ is the $2\times 2$ matrix operating on ${\mathfrak h}'=\langle x,y\rangle$. Notice that, up to rescaling the parameter $t$, we can assume that either $\beta = 0$ or $\beta = 1$.
Moreover, up to interchanging the basis vector of each $T$-invariant plane in ${\mathfrak h}$, we can assume that $\alpha_i$ is non-negative, for all $i=1,\dots , m$, and, up to interchanging the ordering  of the planes in the basis, we can assume that $\alpha_i\leq \alpha_{i+1}$, for all $i=1,\dots , m-1$.

If ${\mathfrak h}$ contains two proper ideals ${\mathfrak i}_1$ and ${\mathfrak i}_2$ which are invariant under $T$ such that
$[{\mathfrak i}_1, {\mathfrak i}_2]=0$ and ${\mathfrak i}_1 \cap {\mathfrak i}_2={\mathfrak h}'$, we say that ${\mathfrak h}$ is \emph{$T$-decomposable} into the direct sum of ${\mathfrak i}_1$ and ${\mathfrak i}_2$ with amalgamated centre, and we restrict our interest on $T$-undecomposable Lie algebras $\mathfrak h$ of type $\{n,2\}$.
Namely, if $\mathfrak h_1$ and $\mathfrak h_2$ are two $T$-undecomposable nilpotent Lie algebras of type $\{n,2\}$ such that the action of the group
$T$ on the centre $\mathfrak z_1$ of $\mathfrak h_1$ and on $\mathfrak z_2$ of $\mathfrak h_2$ coincides, then the direct sum of $\mathfrak h_1$ and $\mathfrak h_2$ with amalgamated centre is a $T$-decomposable nilpotent Lie algebra  of type $\{n,2\}$, and any $T$-decomposable nilpotent Lie algebra  of type $\{n,2\}$ is obtained in this way.

As above, we define the pair of alternating matrices  $\Big(T_1=(a_{ij}),T_2=(b_{ij})\Big)$ by putting
$$[e_i,e_j]=a_{ij}x+b_{ij}y.$$
Clearly, the Lie algebra $\mathfrak h$ is $T$-decomposable if and only if the matrices $T_1$ and $T_2$ can be put into the same diagonal block form and $T$ leaves invariant the subspaces corresponding to the blocks.

Writing
\begin{equation}\label{partialeven} \begin{array}{c}
\left\{\partial_0(t):= \Big(\alpha_1 t\cdot J \oplus\cdots\oplus\alpha_m t\cdot J\Big)| \; t \in \mathbb R \right\}, \ \hbox{ for $n$ even }  \\
 \left\{\partial_0(t):= \Big(0 \oplus \alpha_1 t\cdot J \oplus\cdots\oplus\alpha_m t\cdot J\Big) |\; t \in \mathbb R \right\}, \hbox{ for $n$ odd, } \end{array} \nonumber \end{equation}
since ${\mathfrak t}$ operates as an algebra of derivations of ${\mathfrak h}$, that is, $$[e_i,e_{j}]^{\partial(t)}=[e_i^{\partial(t)},e_{j}]+[e_i,e_{j}^{\partial(t)}],$$ for a generator of ${\mathfrak t}$ , e. g.  for $t=1$,
 we get
\begin{equation}\label{T1T2partial}\begin{array}{l}
\beta   T_2=\partial_0(1)'T_1+T_1\partial_0(1)\\
-\beta  T_1=\partial_0(1)'T_2+T_2\partial_0(1), \end{array}  \end{equation}
hence we find
\begin{equation}\label{T1T2partialquadratic}\begin{array}{l}
\beta^2  T_1=-\partial_0(1)'^2T_1-2\partial_0(1)'T_1\partial_0(1)-T_1\partial_0(1)^2\\
\beta^2  T_2=-\partial_0(1)'^2T_2-2\partial_0(1)'T_2\partial_0(1)-T_2\partial_0(1)^2. \end{array}  \end{equation}
We arrange the matrices $T_1$ and $T_2$ into $2\times 2$ blocks $A_{hk}$ and $B_{hk}$ with $h, k=1,\dots ,m$ (in the case where $n=2m+1$, we denote the $1\times 2$ blocks of the first row with $A_{0k}$ and $B_{0k}$ and we put $\alpha_0=0$).
Then (\ref{T1T2partial}) is equivalent to
\begin{equation}\label{T1T2}\begin{array}{l}
\beta   B_{hk}=-\alpha_h  J A_{hk}+\alpha_k  A_{hk}J \\
-\beta   A_{hk}=-\alpha_h  J B_{hk}+\alpha_k  B_{hk}J
\end{array}  \end{equation}
and (\ref{T1T2partialquadratic}) is equivalent to
\begin{equation}\label{T1T2blocks}
\begin{array}{l}
-\beta^2\cdot A_{hk}=-(\alpha_h^2+\alpha_k^2)A_{hk}+2\alpha_h\alpha_kJ'A_{hk}J\\
-\beta^2\cdot B_{hk}=-(\alpha_h^2+\alpha_k^2)B_{hk}+2\alpha_h\alpha_kJ'B_{hk}J.
\end{array}
\end{equation}
Notice that, since $\alpha_0=0$, the above equations still hold, with a slight abuse of notation, in the case where $h=0$.

Considering $T_1$ and $T_2$ as split-quaternion matrices, we write $A_{hk}=z_1+z_2\omega$ for suitable complex numbers $z_1$ and $z_2$. Then the equations (\ref{T1T2blocks}) give
\begin{equation}\label{T1T2blocksbeta}\begin{array}{l}
(\alpha_h-\alpha_k)^2z_1=\beta^2z_1, \\
(\alpha_h+\alpha_k)^2z_2=\beta^2z_2. \end{array}  \end{equation}

\begin{Remark}\label{undecomposable}
Notice that, if $h$ is such that, for any $k$ with $\alpha_h\ne\alpha_k$, all the blocks $A_{hk}$ and $B_{hk}$ are zero, then $\mathfrak h$ is $T$-decomposable. By the way, from equations (\ref{T1T2}) it follows that, in the case where $\beta\ne 0$, $A_{hk}$ is zero if and only if $B_{hk}$ is zero.
\end{Remark}

\bigskip
\noindent
\centerline{\bf 3. The case where $\beta =0$} \label{beta=0}

\bigskip
\noindent
This is the case where $T$ is contained in a non-commutative compact group of automorphisms of $\mathfrak h$ (see Remark \ref{compact}).

\begin{Theorem}\label{Thmbeta=0}
With the above notations, let ${\mathfrak h}$ be a $T$-undecomposable Lie algebra of type $\{n,2\}$ and let $\beta =0$. Then $n$ is even and, up to a change of basis, the pair $(T_1,T_2)$ is the real form of the pair of complex matrices  $\left(L_{\frac{n+2}{4}}, R_{\frac{n+2}{4}}\right)^\nabla$, if $n\equiv 2\mbox{ mod }4$, or the real form of the pair of complex matrices  $\left(I_{\frac{n}{4}},J_{\frac{n}{4}}(q)\right)^\nabla$, if $n\equiv 0\mbox{ mod }4$. The group $T$ operates, with respect to the chosen basis, as the group of automorphisms $\exp\big(\partial(t)\big)$, where
$$\partial(t)=t\cdot\Big((\oplus  \alpha J)\Big)\oplus 0.
$$\end{Theorem}
\Pro
From the equations (\ref{T1T2}) we deduce that, if $\alpha_h=0\ne\alpha_k$, then $A_{hk}$ and $B_{hk}$ are zero. As ${\mathfrak h}$ is $T$-undecomposable, it follows from Remark \ref{undecomposable} that $\alpha_h$ is positive for any $h=1,\dots, m$ and $n$ is even. Write $A_{hk}=z_1+z_2\omega$ for suitable complex numbers $z_1$ and $z_2$. Then the equations (\ref{T1T2blocks}) give
$$\begin{array}{l}
(\alpha_h-\alpha_k)^2z_1=0, \\
(\alpha_h+\alpha_k)^2z_2=0. \end{array}  $$
This latter forces $z_2=0$, that is, $A_{hk}$ and $B_{hk}$ are the real form of two complex numbers. Moreover, from the former we obtain that either $A_{hk}$ and $B_{hk}$ are zero, or $\alpha_h=\alpha_k$. As ${\mathfrak h}$ is $T$-undecomposable, we exclude the first case, hence we have that $\partial(t)=t\cdot\Big((\oplus  \alpha J)\Big)\oplus 0$. Since $T_1$ and $T_2$ are the real form of complex $m\times m$ matrices $\widehat{T_1}$ and $\widehat{T_2}$ and $T$ operates on them as the complex scalar matrix $\alpha iI_m$, up to a change of basis in the $m$-dimensional complex space, which leaves $T$ invariant, we can assume that $(\widehat{T_1},\widehat{T_2})$ is either $\left(L_{\frac{n+2}{4}}, R_{\frac{n+2}{4}}\right)^\nabla$ (if $n\equiv 2\mbox{ mod }4$), or $\left(I_{\frac{n}{4}},J_{\frac{n}{4}}(q)\right)^\nabla$ (if $n\equiv 0\mbox{ mod }4$). These are $T$-undecomposable over the real numbers, since the only $T$-invariant real planes are $\Pi_1=\langle e_1,e_2\rangle, \dots ,\Pi_{\frac{n}{2}}=\langle e_{n-1},e_n\rangle$. \qed

\begin{Remark}
Up to rescaling the parameter $t$, we can assume $\alpha = 1$ in the above theorem, but we prefer to leave it, because, in the case where ${\mathfrak h}$ is $T$-decomposable, different values of $\alpha$ can occur.\end{Remark}

\begin{Remark} \label{rmkex}
For $n=4$ and $q=i$, that is, $\big(I_{\frac{n}{4}},J_{\frac{n}{4}}(q)\big)=\big((1),(i)\big)$, we obtain Example \ref{example} in the Introduction. \end{Remark}

\begin{Remark}\label{compact}
Let  ${\mathfrak k}$ be a simple compact algebra of derivations of the nilpotent Lie algebra ${\mathfrak h}$ of type $\{n,2\}$, hence it induces on the $2$-dimensional commutator subalgebra ${\mathfrak h}'$ the null map. Any element in ${\mathfrak k}$ generates a $1$-dimensional compact subalgebra of derivations of ${\mathfrak h}$, thus ${\mathfrak h}$ has the structure given in  Theorem \ref{Thmbeta=0}, and its algebra of derivations can be directly computed.\end{Remark}

\bigskip
\noindent
\centerline{\bf 4. The cases where $\beta\ne 0$} \label{beta ne 0}

\bigskip
\noindent
Up to rescaling the parameter $t$, if $\beta\ne 0$, then we can assume that $\beta =1$. From now on, we need to distinguish the cases where the smallest coefficient $\alpha_h$ is zero or, respectively, smaller, equal or greater than $1/2$. The arguments are more or less the ones we give in the following theorem.

\begin{Theorem}\label{alpha=0}
With the notations given in (\ref{Teven}), if $\beta= 1$ and if the smallest coefficient $\alpha_h$ is zero, then, with respect to a suitable basis of $\mathfrak h$, the group  $T$ operates as the group of automorphisms $\exp\big(\partial(t)\big)$, where \begin{equation}\label{partial-alpha=0}\partial(t)=t\cdot\Big((\oplus 0)\oplus(\oplus  J)\oplus(\oplus 2 J)\oplus\dots\oplus(\oplus l J)\Big)\oplus t\cdot J\end{equation}
with the diagonal blocks $(\oplus i J)$ of dimension $d_i\times d_i$ (with $d_i$ even for $i>0$),
and the $T$-undecomposable Lie algebra $\mathfrak h$ is described by the pair $(T_1,T_2)$, where $T_2=T_1(\oplus J)$ and
\begin{equation}\label{T_1}T_1=\left(\begin{array}{c|c|c|c|c|c}
0&W_0&&&&\\ \hline
-W_0'&0&W_1&&&\\ \hline
&-W_1'&0&W_2&&\\ \hline
&&-W_2'&\ddots&\ddots&\\ \hline
&&&\ddots&0&W_{l}\\ \hline
&&&&-W_{l}'&0
\end{array}\right).\end{equation}
The blocks $W_i$ have dimension $d_i\times d_{i+1}$ and:
\begin{itemize}
\item[i)] the block $W_0$ has, in the most general case, the echelon form
$$W_0=\left(\begin{array}{c|c|cc|c}
0&\oplus L_1 & 0 & 0 & 0\\ \hline
0&0&I_{2s}&0&\Omega_{2s}\\
0&0&0&I_{2t}&0
\end{array}\right),$$
where $L_1=(1,0)$ and $\Omega_{2s}$ is the real form of the split-quaternion matrix $\omega I_s$, $s\geq 0$, $t\geq 0$,
and the first zero columns, as well as the blocks $L_1$, are not necessarily being,
\item[ii)] for any $i>0$, the block $W_i$ is the real form of a complex matrix $\widehat{W}_i$ that can be reduced to the \emph{almost echelon} form
$\widehat{W}_i=(0|\widetilde{I}_{r_i})$ (and, in particular, $\widehat{W}_l=\widetilde{I}_{r_l}$), where $\widetilde{I}_{r_i}$ is obtained by the complex identity matrix by possibly adding zero rows,
\end{itemize}
and such that no two successive columns of $T_1$ of indices $2j-1, 2j$ are both zero.
\end{Theorem}
\Pro  Recall that we have chosen a basis such that $0\leq \alpha_i\leq\alpha_{i+1}$ and that, by the equations (\ref{T1T2}), we get $$A_{hk}=0\iff B_{hk}=0.$$
If the coefficient $\alpha_h$ is zero, then the equations (\ref{T1T2blocks}) become
$$\begin{array}{l}
- A_{hk}=-\alpha_k^2A_{hk}\\
-B_{hk}=-\alpha_k^2B_{hk},
\end{array}
$$
and we see that, if $\alpha_k\ne 1$, then $A_{hk}=B_{hk}=0$. As ${\mathfrak h}$ is $T$-undecomposable, by Remark \ref{undecomposable} we find that the smallest non-zero coefficients must be equal to $1$. For the same reason, the next possible coefficients must be equal to $2$ and so on, that is, the derivations must be of the form given in (\ref{partial-alpha=0}).

Moreover, equations (\ref{T1T2blocksbeta}), which with $\beta =1$ become
$$\begin{array}{l}
(\alpha_h-\alpha_k)^2z_1=z_1, \\
(\alpha_h+\alpha_k)^2z_2=z_2, \end{array}
$$
show that $A_{hk}\ne 0$  only if $|\alpha_h-\alpha_k|= 1$, and that $A_{hk}$ is the real form of a non-zero complex number, as soon as $\; 0<\alpha_h=\alpha_k-1$. Thus $T_1$ must be of the form given in (\ref{T_1})
and the blocks $W_i$  are, for $i>0$, the real form of  complex matrices $\widehat{W}_i$.
By the first of equations (\ref{T1T2}), which for $\beta=1$ gives
$$B_{hk}=-\alpha_h \cdot J A_{hk}+\alpha_k \cdot A_{hk}J,$$
the block of $T_2$, corresponding to $\alpha_h=0$ and $\alpha_k=1$, is equal to $W_0(\oplus J)$, and the same holds for
the blocks corresponding to  $0<\alpha_h=\alpha_k-1$, because in these cases $A_{hk}$ commutes with $J$, hence they are of the form ${W}_i(\oplus J)$ also for $i>0$.

\medskip
Consider a basis change block diagonal matrix of the form $$X=(X_0\oplus X_1\oplus\dots\oplus X_l)\oplus I_2,$$ where, for $i>0$, the blocks $X_i$ are $d_i\times d_i$ matrices that are the real form of complex matrices $\widehat{X}_i$, and notice that it leaves the derivations invariant, and changes the block $W_i$ of $T_1$ into the block
$X_i'W_iX_{i+1}$ (and the corresponding block of $T_2$, accordingly).

In order to reduce the blocks $W_i$ into the form given in the claim, we now perform the following algorithm:
\begin{itemize}
\item[i)] Starting from the last block $W_l$ and working upward one by one, by left multiplication $X_{i}'W_i$ with a suitable matrix $X_{i}$, we can assume that the blocks $W_i$ are reduced to   \emph{lower} echelon form, that is, such that, for $h<k$, the pivot of the $h$-th row is on the \emph{right} of the pivot of the $k$-th row.
Moreover we annihilate the entries also below any pivot (as well as above it). For $i=0$, the matrices $X_0$ and $W_0$ are not necessarily the real forms of complex matrices. For $i>0$, on the contrary, they are.
\item[ii)] Let $i>0$. In order to reduce to zero all the row entries which are on the right of any pivot, we operate on the real form $W_{i}$ of the complex matrix $\widehat{W}_{i}$ by adding to a given (complex) column a linear combination of the  \emph{previous} (complex) columns of $W_{i}$, that is,  by the multiplication
$$\left(\begin{array}{c}W_{i}\\ \hline 0\\ \hline -W_{i+1}'\end{array}\right)X_{i+1}$$
with a suitable matrix $X_{i+1}$ which is the real form of an upper triangular complex matrix with any diagonal (complex) entry equal to $I$ and this operation did not change the lower echelon form of $W_{i+1}$. Notice, in fact, that the transpose of a lower echelon matrix is still a lower echelon matrix. Thus we can assume that
the columns of all blocks $W_i$ are either zero or vectors from the canonical basis (taken in the reverse order).
\item[iii)] We have to distinguish the case where $i=0$, because $W_0$ is not necessarily the real form of a complex matrix, but $X_1$ is such, hence the real matrix $X_1$ operates on \emph{pairs} of columns, with indices $2j-1, 2j$. As $W_0$ is a lower echelon matrix with zeros above and below any pivot, considering the pivot of a row and its position with respect to the pivot of the previous row,  we have the following three cases:
$$\left(\begin{array}{c|cc|c}
\cdots&0&1&\cdots\\
\cdots&1&0&\cdots
\end{array}\right),\;\left(\begin{array}{c|cc|c}
\cdots&0&0&\cdots\\
\cdots&1&a&\cdots
\end{array}\right),\;
\left(\begin{array}{c|cc|c}
\cdots&0&0&\cdots\\
\cdots&0&1&\cdots
\end{array}\right).
$$
As any row is virtually the second row of a two-rows real form of a single complex row, in the last two cases, the second row can be reduced to a vector $e_{2i+1}$ of the real canonical basis,
by multiplying on the right with the real form $X_1$ of an upper triangular complex matrix $\widehat{X}_1$ which has the identity in any (complex) diagonal entry but the one corresponding to the pivot, which has to be \begin{equation}\label{pivot}{\frac{1}{1+a^2}}\left(\begin{array}{cc}
1&-a\\
a&1
\end{array}\right)\quad\mbox{or}\quad\left(\begin{array}{cc}
0&-1\\
1&0
\end{array}\right),\end{equation} respectively. In order to show that we can assume that also the entries below $(1,0)$ are zero, we consider, for instance, the minimal case of the matrix
$$\left(\begin{array}{cc|cc}
0&0&1&0\\ \hline
0&1&0&a\\
1&0&0&b
\end{array}\right).$$
The following multiplications, first on the left
$$\left(\begin{array}{c|cc}
1&0&0\\ \hline
-b&1&0\\
a&0&1
\end{array}\right)\left(\begin{array}{cc|cc}
0&0&1&0\\ \hline
0&1&0&a\\
1&0&0&b
\end{array}\right)=\left(\begin{array}{cc|cc}
0&0&1&0\\ \hline
0&1&-b&a\\
1&0&a&b
\end{array}\right)$$
with a real matrix, and second on the right
$$\left(\begin{array}{cc|cc}
0&0&1&0\\ \hline
0&1&-b&a\\
1&0&a&b
\end{array}\right)\left(\begin{array}{cc|cc}
1&0&-a&-b\\
0&1&b&-a\\ \hline
0&0&1&0\\
0&0&0&1
\end{array}\right)=\left(\begin{array}{cc|cc}
0&0&1&0\\ \hline
0&1&0&0\\
1&0&0&0
\end{array}\right)$$
with the real form of a complex matrix,
show therefore that the block $L_1=(1,0)$ has only zero blocks on left and right, and above and below.

We are left with the first case
$$
\left(\begin{array}{c|cc|c}
\cdots&0&1&\cdots\\
\cdots&1&0&\cdots
\end{array}\right),$$
which is the case where the two rows can be seen as a row with entries in the algebra ${\mathbb H}_-$ of split-quaternions, and with pivot $\omega$. By multiplying on the right with the real form  of an upper triangular complex matrix which has the identity in any (complex) diagonal entry, we reduce each non zero entry $z_1+z_2\omega$ to $z_1$ and, by another multiplication, we can assume that, in the row, only the most leftward entry $z_1$ is non-zero. Finally, we reduce it to $1$ by multiplying on the right with the real form of an upper triangular complex matrix which has $z_1^{-1}$ in the suitable (complex) diagonal entry.

Thus, also in the case $i=0$, we can assume that
the columns of the block $W_0$ are either zero or vectors from the canonical basis.
Moreover this operation did not change the lower echelon form of $W_1$, however some row has been multiplied by a complex scalar at point (\ref{pivot}) and in the above reduction of $z_1$ to $1$. These rows can be reduced again to vectors of the canonical basis, by multiplying on the right with a suitable complex diagonal matrix, and so on with the successive blocks $W_i$.
\item[iv)] Starting again from the last block and working upward, we change now the lower echelon blocks into upper echelon blocks by multiplying on the left with the suitable permutation matrix. We still indicate the blocks by $W_k$.
\item[v)] Starting now from the first block $W_0$, by multiplying on the left by a real permutation matrix $Q_0'$ and on the right by a complex permutation matrix $Q_1$, in the most general case we reduce it to the form
$$Q_0'W_0Q_1=\left(\begin{array}{c|c|cc|c}
0&\oplus L_1 & 0 & 0 & 0\\ \hline
0&0&I_{2s}&0&\Omega_{2s}\\
0&0&0&I_{2t}&0
\end{array}\right),$$
where $\Omega_{2s}$ is the real form of the split-quaternion matrix $\omega I_s$, $s\geq 0$, $t\geq 0$,
and the first zero columns, as well as the blocks $L_1$, are not necessarily being
(notice that $W_0$ cannot have zero rows).
The second block is now $Q_1'W_1$ and its (complex) columns are vectors from the canonical basis, in the order permuted by the multiplication by $Q_1'$. The zero rows are not necessarily at the bottom now, and we cannot move them to the bottom without permuting the columns of $W_0$. On the contrary, by multiplying  on the right by a complex permutation matrix $Q_2$, we can permute the columns and reduce $W_1$ to the \emph{almost echelon} form $\widehat{W}_1=(0|\widetilde{I}_{r_1})$, where $\widetilde{I}_{r_1}$ is obtained by the complex identity matrix by adding zero rows. Going down,  we cannot move the rows of $Q_2'W_2$ without permuting the columns of $W_1$, but we can permute the columns and reduce $\widehat{W}_2$ to the form $\widehat{W}_2=(0|\widetilde{I}_{r_2})$, where $\widetilde{I}_{r_2}$ is obtained by the complex identity matrix by adding zero rows.\end{itemize}
The claim follows after repeating the same argument till the last block $\widehat{W}_l$, which, in particular, will be reduced to the form $\widetilde{I}_{r_l}$, because it cannot have zero (complex) columns. \qed

\begin{Theorem}
With the notations given in (\ref{Teven}), if $\beta=1$ and if the smallest coefficient $\alpha_h$ is greater than 
${\frac{1}{2}}$, then, with respect to a suitable basis of $\mathfrak h$, the group  $T$ operates as the group of automorphisms $\exp\big(\partial(t)\big)$, where
\begin{equation}\label{partial_alpha>1/2beta}\partial(t)=t\cdot\Big((\oplus \alpha J)\oplus(\oplus (\alpha+1) J)\oplus\dots\oplus(\oplus (\alpha+l) J)\Big)\oplus t\cdot J, \nonumber \end{equation}
with the diagonal blocks $(\oplus (\alpha+i-1) J)$ of dimension $d_i\times d_i$ (with $d_i$ even),
then the $T$-undecomposable Lie algebra $\mathfrak h$ is described by the pair $(T_1,T_2)$, where $T_2=T_1(\oplus J)$ and
\begin{equation}T_1=\left(\begin{array}{c|c|c|c|c|c}
0&W_1&&&&\\ \hline
-W_1'&0&W_2&&&\\ \hline
&-W_2'&0&W_3&&\\ \hline
&&-W_3'&\ddots&\ddots&\\ \hline
&&&\ddots&0&W_{l}\\ \hline
&&&&-W_{l}'&0
\end{array}\right).\nonumber \end{equation}
The block $W_i$ has dimension $d_i\times d_{i+1}$ and
is the real form of a complex matrix $\widehat{W}_i$ that can be reduced to the \emph{almost echelon} form
$\widehat{W}_i=(0|\widetilde{I}_{r_i})$ (and, in particular, $\widehat{W}_1=(0|{I}_{r_1})$ and $\widehat{W}_l=\widetilde{I}_{r_l}$), where $\widetilde{I}_{r_i}$ is obtained by the complex identity matrix by adding zero rows
and such that no two successive columns of $T_1$ of indices $2j-1, 2j$ are both zero.
\end{Theorem}
\Pro
Let the smallest coefficient $\alpha_h$ be greater than ${\frac{1}{2}}$ and let $\alpha_k\ne \alpha_h$. By equations (\ref{T1T2blocksbeta}), if $1-\alpha_h\ne\alpha_k\ne 1+\alpha_h$, then $A_{hk}=B_{hk}=0$. Since $\mathfrak h$ is $T$-undecomposable, we have that the closest coefficients are $\alpha_k=1-\alpha_h$ or $\alpha_k=1+\alpha_h$.
But, since $\alpha_h> {\frac{1}{2}}$, we have that $1-\alpha_h<\alpha_h$, thus $\alpha_k=1-\alpha_h$ would be a contradiction to the minimality of $\alpha_h$. Therefore, in this case, the coefficients are $\alpha_h$, $1+\alpha_h$, $2+\alpha_h$, and so on.
With the same arguments as in Theorem \ref{alpha=0}, the claim follows. In this case, also the first block $W_1$ is the real form of a complex matrix $\widehat{W}_1$. \qed

\medskip 
The following case is, somehow, exceptional. In fact, in addition to the various choices of $\widetilde{I}_r$, here it is possible that the real form of an arbitrary skew-symmetric complex matrix gives a class of $T$-undecomposable  Lie algebras $\mathfrak h$.

\begin{Theorem}\label{skewHermitian}
With the notations given in (\ref{Teven}), if $\beta=1$ and if the smallest coefficient $\alpha_h$ is equal to 
${\frac{1}{2}}$, then, with respect to a suitable basis of $\mathfrak h$, the group  $T$ operates as the group of automorphisms $\exp\big(\partial(t)\big)$, where 
\begin{equation}\label{partial_alpha=1/2beta}\partial(t)=t\cdot\Big((\oplus {{\frac{1}{2}}} J)\oplus(\oplus {\frac{3}{2}} J)\oplus\dots\oplus(\oplus {\frac{2l+1}{2}} J)\Big)\oplus t\cdot J.  \nonumber \end{equation}
If we denote by $d_i\times d_i$ ($d_i$ even) the dimension of the block $(\oplus {\frac{2i-1}{2}} J)$,
then the $T$-undecomposable Lie algebra $\mathfrak h$ is described by the pair $(T_1,T_2)$, with $T_2=T_1(\oplus J)$ and
\begin{equation}\label{T_1 1/2} T_1=\left(\begin{array}{c|c|c|c|c|c}
\Omega_{d_1}H&W_1&&&&\\ \hline
-W_1'&0&W_2&&&\\ \hline
&-W_2'&0&W_3&&\\ \hline
&&-W_3'&\ddots&\ddots&\\ \hline
&&&\ddots&0&W_{l}\\ \hline
&&&&-W_{l}'&0
\end{array}\right)\end{equation}
where $\Omega_{d_1}$ is the real form of the split-quaternion matrix $\omega I_{\frac{d_1}{2}}$ and:
\begin{itemize}
\item[i)] $H$ and $W_1$ are the real form of the complex matrices
$$\widehat{H}=\left(\begin{array}{c|cc} \widehat{H}_0&0&0\\ \hline 0&0&I_r\\0&-I_r&0\end{array}\right),\quad\widehat{W}_1=\left(\begin{array}{c|cc} 0&0&I_s\\ \hline 0&0&0\\0&0&0\end{array}\right),$$
respectively, with $s\geq 0$, $r\geq 0$, $d_1=2(s+2r)$, and $\widehat{H}_0$ is a $s\times s$ skew-symmetric complex matrix,
\item[ii)] for any $i>1$, the block $W_i$ of dimension $d_i\times d_{i+1}$  is the real form of a complex matrix $\widehat{W}_i$ that can be reduced to the \emph{almost echelon} form
$\widehat{W}_i=(0|\widetilde{I}_{r_i})$ (and, in particular, $\widehat{W}_l=\widetilde{I}_{r_l}$), where $\widetilde{I}_{r_k}$ is obtained by the complex identity matrix by adding zero rows,
\end{itemize}
and such that no two successive columns of $T_1$ of indices $2j-1, 2j$ are both zero.
\end{Theorem}
\Pro
If $\alpha_h= {\frac{1}{2}}$, then the assumption $\alpha_k=1-\alpha_h$ leads to a contradiction to $\alpha_k\ne\alpha_h$. Thus, also in this case, the distinct coefficients are $\alpha_h$, $1+\alpha_h$, $2+\alpha_h$, and so on. By equations (\ref{T1T2blocksbeta}), we see, however, that, if $\alpha_k=\alpha_h$, then $A_{hk}=z_2\omega$, hence the $d_1\times d_1$ block corresponding to the values of $\alpha_k$ equal to ${\frac{1}{2}}$ is the real form of a split-quaternion matrix $\omega \widehat{H}$. Notice that, whereas the real form of $\widehat{H}$  is skew-symmetric if and only if  $\widehat{H}$ is skew-Hermitian, the real form of $\omega \widehat{H}$ is skew-symmetric if and only if  $\widehat{H}$ is skew-symmetric. 
This forces $T_1$ and $T_2$ to be of the form given in (\ref{T_1 1/2}).

With the same argument of Theorem \ref{alpha=0}, we reduce $W_i$ to the almost echelon form $(0|\widehat{I}_{r_i})$,
and in particular $W_l=\widehat{I}_{r_l}$ and, since we can operate on the first row by multiplication on the left with a further matrix, we reduce $W_1$ to the real form of
$${\widehat{W}}_1=\left(\begin{array}{c|c} 0&I_{s}\\ \hline 0&0\end{array}\right).$$
A basis change matrix of the form $X_1\oplus I_{d_2}\oplus\dots\oplus I_{d_l}$ transforms the blocks $W_1$ and $\Omega_{d_1}H$ into $X_1'W_1$ and $X_1'\Omega_{d_1}HX_1$, respectively. 

\noindent 
In order to leave invariant the echelon form of $W_1$, we have to take
$\widehat{X}_1=\left(\begin{array}{c|c} I_s&0\\ \hline C&D\end{array}\right)$.

Notice now that $X_1'\Omega_{d_1}HX_1$ is the real form of
$$\widehat{X}_1^\dag \omega I_{\frac{d_1}{2}}\widehat{H}\widehat{X}_1= \omega I_{\frac{d_1}{2}}\widehat{X}_1'\widehat{H}\widehat{X}_1.$$
If we write $\widehat{H}=\left(\begin{array}{c|c} H_0&H_1\\ \hline -H_1'&H_2\end{array}\right)$, with $H_0$ of dimension $s\times s$, we see that the congruence $\widehat{X}_1'\widehat{H}\widehat{X}_1$ changes $H_2$ into $D' H_2D$ and $H_1$ into $(H_1+C' H_2)D$. 

Thus we can reduce, firstly, $H_2$ to the canonical form $\left(\begin{array}{c|c} 0&I_r\\ \hline -I_r&0\end{array}\right)$, because $H_2$ has to be non-degenerate, or $T_1$ would have a zero (complex) row. Finally, since $H_2$ is non-degenerate, we reduce $H_1$ to zero, taking  $C'=-H_1H_2^{-1}$.
\qed

\begin{Remark} In the case where $\partial(t)=(\oplus {{\frac{1}{2}}}t\cdot J)\oplus t\cdot J$, we find
that $(T_1,T_2)=(\oplus I_2,\oplus J)^\nabla$.
\end{Remark}
In the following last theorem we will change the ordering of the coefficients $\alpha_h$ defining the derivations $\partial(t)$.

\begin{Theorem}\label{last}
With the notations given in (\ref{Teven}), if $\beta=1$ and if the smallest coefficient $\alpha_h$ is smaller than ${\frac{1}{2}}$, then, with respect to a suitable basis of $\mathfrak h$, the group  $T$ operates, in the most general case, as the group of automorphisms $\exp\big(\partial(t)\big)$, where $\partial(t)=\partial_1(t)\oplus\partial_2(t)\oplus t\cdot J$ with \begin{equation}\label{partial1}\partial_1(t)=t\cdot\Big((\oplus (l_1-\alpha) J)\oplus\dots\oplus(\oplus(2-\alpha)J)\oplus(\oplus(1-\alpha)J)\Big) \nonumber \end{equation}
\begin{equation}\label{partial2}\partial_2(t)=t\cdot\Big((\oplus\alpha J)\oplus(\oplus (\alpha+1) J)\oplus\dots\oplus(\oplus (\alpha+l_2) J)\Big) \nonumber \end{equation}
and the $T$-undecomposable Lie algebra $\mathfrak h$ is described by the pair $(T_1,T_2)$, with $T_2=T_1(\oplus J)$ and
\begin{equation}\label{T_1 <1/2}T_1=\left(\begin{array}{c|c|c|c|c|c|c|c}
0&V_{l_1-1}&&&&&&\\ \hline
-V_{l_1-1}'&\ddots&\ddots&&&&&\\ \hline
&\ddots&0&V_1&&&&\\ \hline
&&-V_1'&0&\Omega_{d_1}S&&&\\ \hline
&&&-S'\Omega_{d_1}&0&W_{1}&&\\ \hline
&&&&-W_{1}'&0&\ddots&\\ \hline
&&&&&\ddots&\ddots&W_{l_2}\\ \hline
&&&&&&-W_{l_2}'&0
\end{array}\right)\end{equation}
where the blocks $V_i$, $W_i$ and $S$ are the real form of  complex matrices $\widehat{V}_i$, $\widehat{W}_i$ and $\widehat{S}$ that can be reduced to the \emph{almost echelon} form
$\widehat{V}_i=(0|\widetilde{I}_{r_i})$, $\widehat{W}_i=(0|\widetilde{I}_{s_i})$ and $\widehat{S}=(0|\widetilde{I}_{s})$  (and, in particular, $\widehat{V}_{l_1-1}=(0|{I}_{r_{l_1-1}})$ and $\widehat{W}_{l_2}=\widetilde{I}_{r_{l_2}}$), where $\widetilde{I}_{r_k}$ is obtained by the complex identity matrix $I_l$ by adding zero rows, $\Omega_{d_1}$ is the real form of the split-quaternion matrix $\omega I_{\frac{d_1}{2}}$,
and such that no two successive columns of $T_1$ of indices $2j-1, 2j$ are both zero.
\end{Theorem}
\Pro
At last, let the smallest coefficient $\alpha_h$ be smaller than ${\frac{1}{2}}$, thus the closest possible coefficients are $\alpha_k=1-\alpha_h$
and $\alpha_k=1+\alpha_h$. The latter gives, as above, the coefficients $\alpha_h$, $1+\alpha_h$, $2+\alpha_h$, and so on. The former gives moreover the coefficients $1-\alpha_h$, $2-\alpha_h$, and so on, and no other, since $(1-\alpha_h)-1$ is negative and $1-(1-\alpha_h)$ is again $\alpha_h$.
Thus, in this case, the coefficients are
$$l_1-\alpha,\dots, 2-\alpha,1-\alpha,\alpha,1+\alpha,2+\alpha,\dots,l_2+\alpha.$$
Notice that no coefficient of the form $a-\alpha_h$ can be equal to $b+\alpha_h$, because $\alpha_h< {\frac{1}{2}}$.
It follows that:
\begin{itemize}
\item[i)] if $\alpha_h=a+\alpha$ and $\alpha_k=b+\alpha$, then $\alpha_h+\alpha_k$ is not an integer and, by equations (\ref{T1T2blocksbeta}), we see that $A_{hk}$ is non-zero only if $|b-a|=1$ and that $A_{hk}$ is the real form of a complex number,
\item[ii)] if $\alpha_h=a+\alpha$ and $\alpha_k=b-\alpha$, then $|\alpha_h-\alpha_k|$ is not an integer and, by equations (\ref{T1T2blocksbeta}), we see that $A_{hk}$ is non-zero only if $a+b=1$, that is $\alpha_h=\alpha$ and $\alpha_k=1-\alpha$, and that  $A_{hk}$ is the real form of a split-quaternion $\omega z_2$,
\item[iii)] if $\alpha_h=a-\alpha$ and $\alpha_k=b-\alpha$, then $\alpha_h+\alpha_k$ is not an integer and, by equations (\ref{T1T2blocksbeta}), we see that again $A_{hk}$ is non-zero only if $|b-a|=1$ and that $A_{hk}$ is the real form of a complex number.
\end{itemize}
Thus, $T_1$ is of the form given in (\ref{T_1 <1/2}). With the same arguments as in  Theorem \ref{alpha=0}, we can reduce the blocks $V_i$, $W_i$ and $S$ to the form given in the claim.
\qed

\begin{Remark} In Theorem \ref{last} it can happen that $\partial(t)=\Big(\partial_1(t)\oplus (\oplus\alpha t\cdot J)\Big)\oplus t\cdot J$ and
$$T_1=\left(\begin{array}{c|c|c|c|c}
0&V_{l_1-1}&&&\\ \hline
-V_{l_1-1}'&\ddots&\ddots&&\\ \hline
&\ddots&0&V_1&\\ \hline
&&-V_1'&0&\Omega_{d_1}S\\ \hline
&&&-S'\Omega_{d_1}&0
\end{array}\right)$$
or that $\partial(t)=\partial_2(t)\oplus t\cdot J$
and $$T_1=\left(\begin{array}{c|c|c|c}
0&W_{1}&&\\ \hline
-W_{1}'&0&\ddots&\\ \hline
&\ddots&\ddots&W_{l_2}\\ \hline
&&-W_{l_2}'&0
\end{array}\right).$$
\end{Remark}

\medskip
\noindent
Giovanni Falcone, Dipartimento di Matematica e Informatica, via Archirafi 34, I-90123 Palermo (Italy), giovanni.falcone@unipa.it

\smallskip
\noindent
\'Agota Figula, University of Debrecen, Institute of Mathematics
\newline
\noindent
H-4010 Debrecen, P.O.Box 12, Hungary, figula@science.unideb.hu


\begin{thebibliography}{9999}

\bibitem{Bartolone} C. Bartolone, A. Di Bartolo and G. Falcone, Nilpotent Lie algebras with $2$-dimensional commutator ideals, Linear Algebra Appl. {\bf 434} (2011), 650-656.

\bibitem{Bartolone2} C. Bartolone, A. Di Bartolo and G. Falcone, Derivations of a nilpotent Lie algebra with $2$-dimensional commutator ideal and $1$-dimensional centre, submitted.

\bibitem{Belitskii} G. Belitskii, R. Lipyanski and V. V. Sergeichuk, Problems of classifying associative or Lie algebras and triples of symmetric or skew-symmetric matrices are wild, Linear Algebra Appl. {\bf 407} (2005), 249-262.

\bibitem{Olomuc} G. Falcone and M. A. Vaccaro, Kronecker Modules and Reductions of a
Pair of Bilinear Forms, Acta Univ. Palacki. Olomuc., Fac. rer. nat., Mathematica {\bf 43} (2004), 55-60.

\bibitem{Gauger} M. A. Gauger, On the classification of metabelian Lie algebras, Trans. Amer. Math. Soc. {\bf 179} (1973), 293-329.

\bibitem{Kronecker} L. Kronecker, Algebraische Reduktion der Scharen bilinearer Formen, Sit\-zungs\-ber. Akad. Berlin (1890), 763-776.


\bibitem{Scharlau} R. Scharlau, Paare alternierender Formen, Math. Z., {\bf 147} (1976), 13-19.

\bibitem{tamaru} H. Tamaru and H. Yoshida, Lie groups locally isomorphic to generalized Heisenberg groups, Proc. Amer. Math. Soc.  
{\bf 136} (2008), 3247-3254. 

\bibitem{Vergne} M. Vergne, Cohomologie des alg\'ebres de Lie nilpotentes. Application \'a l'\'etude de la vari\'et\'e des algebres de Lie nilpotentes, Bull. Soc. Math. France {\bf 98} (1970), 81-116.

\bibitem{wilson} E. N. Wilson, Isometry groups on homogeneous nilmanifolds, Geom. Dedicata {\bf 12} (1982), 337-346. 





\end{thebibliography}
\end{document}